\documentclass{article}
\usepackage[utf8]{inputenc}

\usepackage{palatino}

\usepackage{mdframed}

\usepackage{comment}

\usepackage{amsfonts, amssymb}

\usepackage{amsmath}

\usepackage{mdframed}

\usepackage{hyperref}

\usepackage{todonotes}

\usepackage{comment}

\newcommand*{\defeq}{\stackrel{\text{def}}{=}}

\newtheorem{question}{Question}

\title{Higher-Order Platonism and Multiversism}

\author{Claudio Ternullo\footnote{Department of Mathematics and Computer Science, University of Barcelona (Spain). \ttfamily \textbf{claudio.ternullo@ub.edu}\rmfamily. For the writing of this paper, I have been supported by the AGAUR and the European Union through a Beatriu de Pin\'os (Marie-Sk\l odowska Curie COFUND) Fellowship (grant n. 00192 BP 2018). I am indebted to Edward Zalta for valuable comments and further clarifications of his theory, and two anonymous reviewers for their extensive, and detailed, feedback.}}

\date{}

\begin{document}

\maketitle

\begin{abstract}

\noindent
Joel Hamkins has described his multiverse position as being one of `higher-order realism -- Platonism about universes', whereby one takes \textit{models of set theory} to be actually existing objects (vis-à-vis `first-order realism', which takes only \textit{sets} to be actually existing objects). My goal in this paper is to make sense of the view in the very context of Hamkins' own multiversism. To this end, I will explain what may be considered the central features of higher-order platonism, and then will focus on Zalta and Linsky's Object Theory, which, I will argue, is able to faithfully express Hamkins' conception. I will then show how the embedding of higher-order platonism into Object Theory may help the Hamkinsian multiversist to respond to salient criticisms of the multiverse conception, especially those relating to its articulation, skeptical attitude, and relationship with set-theoretic practice. 
    
\end{abstract}

\section{Hamkins' Multiverse and Platonism}\label{Intro}

The `set-theoretic multiverse conception' holds that the subject matter of set theory is not a \textit{single} universe (in particular, not the well-founded, cumulative hierarchy $V$), but \textit{some} (or \textit{all}) the \textit{models} of some set theory $T$ (most commonly, of $\mathsf{ZFC}$).\footnote{I assume that readers have some familiarity with `models of set theory'. For a systematic treatment, which covers all models mentioned in this paper, see \cite{kunen2011}.} 

The most lucid and radical expression of this conception can be found in \cite{hamkins2012}. This is not to deny that alternative conceptions are possible (and, indeed, have been formulated); only, Hamkins' multiversism seems to stand out as the most coherent, and self-conscious, available form of \textit{set-theoretic pluralism}.\footnote{The literature on `multiverse conceptions' is vast. A rough classification, and discussion of the main features of each conception, may be found in \cite{afht2015}. For an introduction to, and discussion of, the basics of `set-theoretic pluralism', see \cite{linnebo2017}, section 12.4.}  

Despite having the unquestionable merit of bringing about the emergence of a new approach to set-theoretic foundations, \cite{hamkins2012} plainly rests on philosophically controversial claims. The invocation of the existence of a \textit{plethora} of universes of set theory looks especially controversial.

Moreover, as a preliminary characterisation of the  philosophy it aims to advocate throughout, the paper makes clear, from the beginning, that the `multiverse position' is: 

\begin{quote}

..one of higher-order realism—Platonism about universes—and I defend it as a realist position asserting actual existence of the alternative set-theoretic universes into which our mathematical tools have allowed us to glimpse. The multiverse view, therefore, does not reduce via proof to a brand of formalism. In particular, we may prefer some of the universes in the multiverse to others, and there is no obligation to consider them all as somehow equal.\footnote{As pointed out by one reviewer, also formalists may want to pick out some particular universe as their \textit{preferred} universe of sets. An instance of this attitude can be found, for instance, in \cite{shelah2003}, p. 211f.} \\ (\cite{hamkins2012}, p. 417)
    
\end{quote}

\noindent
In a fully consequential manner, a substantial amount of work in the paper is, then, devoted to articulating the idea that \textit{set-theoretic models}, in particular, models obtained through \textit{forcing}, may be construed as `objects' in the \textit{platonic} sense.\footnote{For a technical exposition of \textit{forcing}, cf. again the mentioned \cite{kunen2011}, or \cite{jech2003}, Ch. 14.} 

The goal of this paper is to investigate whether and how Hamkins' multiverse conception may legitimately claim to be able to do that, and on what metaphysical grounds; a parallel question, which I will also investigate, is that of how well `higher-order realism' serves Hamkins' own purpose of laying out a \textit{pluralist} foundation of set theory.   


In order to be clear about exactly what is the position I shall be dealing with, I shall first explicate what can be viewed as the two fundamental tenets of Hamkins' conception, which I shall henceforth indicate as Higher-Order Platonism (HOP). We may summarise them as follows: 

\vspace{11pt}

\noindent
\textbf{Platonism (PLAT)}. The models of set theory are \textit{platonic} entities, i.e., are self-standing, independently existing, physically acausal, mathematical objects. 

\vspace{11pt}

\noindent
\textbf{Perspectivism (PERSP)}. The set-theoretic multiverse is always \textit{relative to} a single model, that is, it is \textit{describable} only from the point of view of one particular set-theoretic model.

\vspace{11pt}

\noindent
(PLAT) already looks very peculiar, but I am inclined to see (PERSP) as being more crucial to the HOP-ist's purposes. Some further explication of (PERSP) may, therefore, already at this stage, be useful. 

The view implies that set-theoretic structures are not definable (`observable') from the point of view of an absolute background structure.\footnote{The reader is warned that `model' and `structure' are used  interchangeably throughout.} For instance, one usually thinks of inner models like $L(\mathbb{R})$, HOD, etc., as \textit{just} subclasses of $V$. But as Hamkins himself explains:

\begin{quote}

I counter this attitude, however, by pointing out that much of our knowledge of these inner models has actually arisen by considering them inside various outer models. We understand the coquettish nature of HOD, for example, by observing it to embrace an entire forcing extension, where sets have been made definable, before relaxing again in a subsequent extension, where they are no longer definable (\cite{hamkins2012}, p. 418)
    
\end{quote}

\noindent
So, it is not simply through viewing them as \textit{definable} subclasses of $V$ that one comes to understand the very nature of inner models, but rather through studying their behaviour within \textit{other} models, such as, for instance, forcing extensions of another (provisional) background universe $V$. 
This is an instance of what I call `perspectivism' with respect to models (with respect to set theory at large): a single model will display different, possibly more interesting, features if studied in the context of (as part of) other models. 

In an effort to further clarify this view, more recently Hamkins has suggested the following characterisation of (PERSP): 

\begin{quote}

The multiverse perspective ultimately provides what I
view as an enlargement of the theory/metatheory
distinction. There are not merely two sides of this
distinction, the object theory and the meta-theory; rather,
there is a vast hierarchy of metatheories. Every set-theoretic
context, after all, provides in effect a
metatheoretic background for the models and theories that
exist in that context -- a model theory for the models and
theories one finds there. \\ (\cite{hamkins2020}, pp. 297-8)
    
\end{quote}

\noindent
An equally promising way of construing (PERSP), then, is through envisioning the existence of a \textit{hierarchy} of set theories, each reflecting a (provisional) theory/metatheory distinction with respect to another one,  such that: the model theory of a theory $T$ (that is, the `multiverse' of $T$) is observed from the point of view of another theory (and provisional \textit{metatheory}) $T'$, the multiverse of $T'$ from the point of view of the theory (and provisional \textit{metametatheory}) $T''$ and so on.\footnote{Note that this hierarchy of theories should not be taken to correspond to the hierarchy of set theories linearly ordered by \textit{consistency strength} as calibrated by large cardinal axioms.}
On pain of contradicting the essence of (PERSP), there can be no \textit{ultimate} level in this hierarchy, insofar as there is no ultimate background theory from whose point of view one can `observe' the model theory of other theories.

Now, (PLAT) and (PERSP) may, to some extent, seem to contradict each other. After all, many `classic' platonists take platonic entities to be \textit{determinate} and \textit{static}.\footnote{But then the issue of referential indeterminacy, as illustrated by \cite{putnam1980}, arises for this kind of platonists. A quick response to Putnam's classic arguments would consist in taking platonic entities to be determinate in an idealised, pre-theoretic sense. Another one would be to thrive on the \textit{categoricity} of higher-order axioms; for an exhaustive overview of, and troubles with, this strategy, see \cite{button-walsh2018}.}
How is it, then, that they change their `nature' depending on which model, as postulated by (PERSP), one observes them from? By (PERSP), for instance, a set-theoretic structure, say, $\mathfrak{A}$, could, if observed from the point of view of another structure $\mathfrak{A'}$, look different from the way it could look if observed from the perspective of one further structure $\mathfrak{A''}$.

One major issue relating to HOP, is, as a consequence, its \textit{coherence}, i.e., whether it is able to yield a conceptual framework which makes sense of the co-existence between (PLAT) and (PERSP). 

What I aim to do in this paper is precisely to describe a metaphysical framework which does not only account for the coherence of HOP, but is also able to provide responses to several concerns about Hamkins' multiverse which have been left unanswered.

\section{Higher-Order Platonism: Main Features}\label{HOP}

Before doing that, I need to discuss a few further features of HOP, in particular, I should make clear in what sense it could be seen as instantiating \textit{mathematical (set-theoretic) platonism}. 

Although accounts of platonism may differ, scholars usually agree on the its basic tenets, which are taken to consist in the following claims:\footnote{For an overview of platonism in mathematics, see \cite{linnebo2018}, which I have used here to characterise the fundamentals of the position.} 

\vspace{11pt}

\noindent
\textbf{Existence (E)}. There \textit{exist} mathematical objects.

\vspace{11pt}

\noindent
\textbf{Abstractness (A)}. Mathematical objects are \textit{abstract} (as opposed to \textit{physically instantiated}, \textit{concrete}) objects.\footnote{Abstractness, construed as \textit{non-spatiotemporality}, may not be a peculiar and exclusive feature of mathematical entities. In particular, also some physical particles, as, for instance, described in quantum field theories, may not be spatio-temporally located. I am indebted to an anonymous reviewer for pointing this fact to me.}

\vspace{11pt}

\noindent
\textbf{Independence (I)}. Mathematical objects are \textit{independent} of minds, language, uses and practices. 

\vspace{11pt}

The E, A, I triad also seems to fit in with HOP: that there exist mathematical objects is implicit in Hamkins' `platonism about universes', only the objects referred to by HOP are not the `classic', first-order mathematical objects; models of set theory seem to be no less \textit{abstract} than first-order objects, and, finally, Independence may carry over to HOP without much ado or, in any case, it does not seem to be less appropriate to HOP than it is to classic (first-order) platonism.        

Some platonists also hold truth-value realism (TVR), the contention that all set-theoretic statements have a \textit{determinate} truth-value. Now, if one takes TVR to be an essential component of classic platonism as much as the E, A, I triad, then HOP noticeably departs from classic platonism, since, as is clear, each set-theoretic model will fix the truth of some set-theoretic statements in a way which might differ from that of another model.\footnote{For noticeable examples of platonists also holding TVR, see \cite{godel1947} and \cite{martin2001}. \cite{isaacson2011}, while subscribing to TVR up to the level of $\mathsf{ZFC}$, holds that at least \textit{some} set-theoretic statements (e.g., the existence of large cardinals), are genuinely \textit{indeterminate}.}  Hence, by HOP, it is not always the case that, for a given statement $\phi$, $\phi$ is \textit{true} or \textit{false}: in many cases $\phi$ will \textit{neither} be true \textit{nor} false (it will be \textit{indeterminate}).  

However, one may attempt to rescue TVR by \textit{relativising} it to a specific structure. One could, that is, propose a revised version of TVR along the following lines: 

\vspace{11pt}

\noindent
\textbf{RTVR (Relativised Truth-Value Realism)}. Each mathematical (set-theoretic) statement is \textit{determinately} true or false \textit{in} a given mathematical structure. 

\vspace{11pt}

\noindent
Moreover, one could also attempt to go beyond the Tarskian concept of `truth in a model', and define a concept of `truth \textit{in all} models (of some kind)'. Such an account of truth is, for instance, described by \cite{woodin2011c} in connection with the `set-generic multiverse conception'. The multiverse concept of truth is as follows:

\vspace{11pt}

\noindent
\textbf{Multiverse Conception of Truth (Falsity)}. $\phi$ is true (false) iff it is true (false) in \textit{all} members of the set-generic multiverse, that is, the collection of all set-generic extensions of an initial model $M$, and of their grounds.\footnote{This notion is best explicated in terms of `satisfaction' in Woodin's own $\Omega$-logic, `$\models_{\Omega}$', for which also see \cite{woodin2001} or \cite{bagaria2005}.} 

\vspace{11pt}

\noindent
However, this strategy does not seem to square well with HOP. This is clear if one considers this further feature of Hamkins' multiverse, as expressed in the following quote:

\begin{quote}
    The background idea of the multiverse, of course, is that there sho-uld be a large collection of universes, each a model of (some kind of) set theory. There seems to be no reason to restrict inclusion only to $\mathsf{ZFC}$ models, as we can include models of weaker theories $\mathsf{ZF}$, $\mathsf{ZF^{-}}$, $\mathsf{KP}$, and so on, perhaps even down to second-order number theory, as this is set-theoretic in a sense (\cite{hamkins2012}, p. 436).
\end{quote}

\noindent
What can be gleaned from the quote above is that Hamkins' multiverse is not exhausted by any collection of set-theoretic models (of any kind). At this point, this would hardly strike anyone as surprising: already by (PERSP), every single universe must `yield' its own multiverse, and it should not be possible to amalgamate all such multiverses into a unique, coherent model-theoretic super-structure. As a consequence, the Hamkinsian multiversist ought to abandon entirely the goal of defining `multiverse truth' in the way sketched above, although, in principle, he might subscribe to RTVR.    

But the quote also provides us with another fundamental piece of information about HOP: all models of \textit{all theories} (of sets) are members of Hamkins' multiverse. This means that the Hamkinsian multiversist is not only a \textit{semantic pluralist}, in the sense that he thinks that truth may vary across the multiverse, but also a \textit{proof-theoretic pluralist}, insofar as he thinks that any collection of consistent set-theoretic axioms provides us with some version of `set-theoretic truth' (that it meets, as Hamkins says himself, some `concept of set').

I will discuss in depth, and actively exploit, the entangled character of the two positions later on; for the time being, it is worth summarising this further fundamental feature of HOP as follows:

\vspace{11pt}

\noindent
\textbf{Proof-Theoretic Pluralism (PTP)}. All \textit{consistent} theories of sets instantiate equal-ly acceptable versions of set-theoretic truth (`concepts of set').

\vspace{11pt}

\noindent
One fundamental complication relating to (PTP) is that, in violation of (PERSP), it could be taken to be an attempt to articulate the whole (or, for that matter, some) multiverse, but then the problem would arise of what the background theory is, from whose viewpoint (PTP) is formulated. This may have awkward consequences. For instance, from the point of view of the theory $\mathsf{ZF}+\neg Con(\mathsf{ZF})$, $\mathsf{ZF}$ has no (standard) models, so, from the point of view of that theory, very oddly, and unexpectedly, the multiverse would not contain any (standard) $\mathsf{ZF}$-models.

One possible response to this issue is provided by \cite{clarke-doane2022}, which proposes to replace `consistent' with `arithmetically sound' in (PTP), a move which would prevent theories like $\mathsf{ZF}+\neg Con(\mathsf{ZF})$ from being seen as acceptable.\footnote{Cf. \cite{clarke-doane2022}, p. 34.}  

In section \ref{SRI}, we shall see that the use of OT can provide us with an alternative and, in my view, equally workable response to this problem. Meanwhile, I will maintain that, at this stage, (PTP) may be seen as enjoying a pre-formalised, informal, status which, \textit{prima facie}, does not rule out inclusion of any theory (and models thereof) in the Hamkinsian multiverse.

\medskip

Our description of HOP is now complete. We have seen that classic platonism's E, A, I triad is kept by HOP, whereas unrelativised semantic determinacy (TVR) is dropped. Moreover, HOP's \textit{semantic pluralism} is also subtly related to \textit{proof-theoretic pluralism},  in a word, the idea that there are as many set-theoretic concepts and versions of truth as, arguably, formal theories of sets.

\section{Object Theory}\label{Forms}

The metaphysical framework which seems most suitable to make sense of HOP is that of \textit{plenitudinous platonism}, which substantially predates the emergence of Hamkins' multiverse conception. Several forms of plenitudinous platonism are available on the market.\footnote{\label{plen}\cite{lewis1986}'s modal realism, \cite{balaguer1995}'s full-blooded platonism,  \cite{shapiro1997}'s \textit{ante rem} structuralism, and \cite{blechschmidt2022}'s topos-theoretic multiversism may all be seen as instances of plenitudinous platonism.} All of these presuppose that, contrary to classic platonism, there exist \textit{plenitudes} of platonic mathematical objects, possibly one plenitude for each item of formal languages: objects, theories, relationships, worlds, models, etc. 

Moreover, the classic platonist is (mostly) construed as believing in a special form of \textit{intuition}, which would provide her with knowledge of mathematical objects as certain (and determinate) as knowledge of physical objects.\footnote{For classic platonism, see \cite{bernays1935}, and \cite{godel1947} (but also see G\"odel's further arguments (and explanations) in \cite{godel1951}). For a useful overview of G\"odel's conception, see \cite{parsons1995}.} Plenitudinous versions of platonism, on the other hand, purport to be able to eliminate the somewhat mysterious character of such an intuition. This is because they hold that what one really has access to when one refers to \textit{abstracta} are totalities of them, which are made epistemically accessible to one through some Comprehension Principle, or, more simply, by the Completeness Theorem for first-order theories.\footnote{For the connections between consistency and plenitudinousness, cf., in particular, \cite{balaguer1995}, pp. 304ff.}

Now, the specific form of plenitudinous platonism with which I will be concerned in this work is Object Theory (OT).\footnote{The material which follows is mostly based on \cite{zalta1983}, \cite{linsky-zalta1995} and \cite{nodelman-zalta2014}. Further details on OT may be found in \cite{zalta2000}, and \cite{linsky-zalta2006}.
For an extensive discussion of OT, see \cite{panza-sereni2013}, section 5.2.} As far as knowledge of mathematical objects (as well as reference to them) is concerned, its authors make clear that this is conveyed by \textit{descriptions}. The authors say:

\begin{quote}

Knowledge of particular abstract objects does not require any
causal connection to them, but we know them on a one-to-one basis
because de re knowledge of abstracta is by description. All one has to do to become so acquainted de re with an abstract object is to understand its descriptive, defining condition, for the properties that an abstract object encodes are precisely those expressed by their defining conditions. So our cognitive faculty for acquiring knowledge of abstracta is simply the one we use to understand the comprehension principle. (\cite{linsky-zalta1995}, p. 547)
 
\end{quote}

\medskip

\noindent
`Encoding', referred to in the quote above, lies at the heart of OT, so it should be explained straight away. This consists in a new form of \textit{predication}, denoted `$xF$', which means: `$x$ encodes $F$' (as opposed to the well-known `$Fx$', which means: `$x$ exemplifies $F$'). The difference between exemplification and encoding can be explained as follows.

An object $x$ exemplifies a property $P$ if that property can be \textit{predicated} of $x$: `2 is prime' is a case in point. An object $x$ encodes a property $P$ if $x$ \textit{is} that property. `2 encodes \textit{primality}' is a case in point. Abstract objects will both exemplify and encode properties.\footnote{A way to characterise `encoding', as opposed to `exemplification', suggested to me by one reviewer, is through viewing encoding as expressing \textit{essences}: so, `$x$ encodes $P$' may be paraphrased as expressing that `$x$ is some \textit{essence} $P$'.}

Encoding is regulated by the following Comprehension Principle:\footnote{OT's language is that of the second-order quantified S5 logic, enriched with the abstraction operator $\lambda$, the unary predicate $E!$, which, applied to the variable $x$, means: `$x$ is concrete' (predicates for `abstract' and `ordinary', resp. $A!$, and $O!$ are, in turn, obtained from $E!$) and the corresponding equality symbol $=_{E!}$. Cf. \cite{nodelman-zalta2014}, pp. 42-3.}

\[(\exists x)(A!x \wedge \forall F(xF \leftrightarrow \phi)) \tag{Comp} \]

\noindent
where $A!x$ means: `$x$ is abstract', $F$ is a property, $\phi$ is a set of conditions (expressing the property $F$). So, (Comp) asserts the existence of a plenitude of abstract objects, each of which encodes all and only those properties expressed by a set of conditions $\phi$. In simpler terms, an abstract object is just an object which \textit{encodes} a collection of properties.\footnote{Alternatively, for any non-empty collection of properties, there is an abstract object which encodes them.} (Comp) is OT's archetypal Plenitude Principle. It allows the formation of \textit{any kind of} abstract objects, in particular, mathematical objects.  

But OT is not exhausted by (Comp). The theory extends \textit{abstractness}, and \textit{plenitudinousness}, to all other `objectifiable' parts of formal languages, among these, crucially, to \textit{theories}.

Formally, a \textit{theory} is just a collection of propositions. In turn, \textit{propositions} are well-formed formulae, consisting of symbols for constants, predicates, relations, all of which are also distinct abstract objects. So, theories, qua abstract objects, are combinations of other abstract objects.  More specifically, through the corresponding Comprehension Principles, one first forms a plenitude of \textit{relations (predicates)}, then a plenitude of \textit{propositions}, construed as 0-place relations. Finally, one forms a plenitude of theories by identifying one single theory $T$ with its \textit{theorems}, that is, with the closure set of the `$\vdash$' relationship, or, \textit{semantically}, with the set of the \textit{logical consequences} of $T$.\footnote{The switch between syntax and semantics is fundamental for Linsky and Zalta's purposes. This can be accounted for, in OT, through using what the authors call Importation Rule, described in \cite{nodelman-zalta2014}, p. 48.} The resulting Comprehension Principle is as follows:

\[ T \defeq \iota x(A!x \wedge (\forall F)(xF \leftrightarrow \exists p (x \models p \wedge F = [\lambda y \ p]))) \tag{Comp$_T$} \] 

\noindent
that is, a theory $T$ is the \textit{only} abstract object which encodes all and only those properties asserted by its \textit{true} propositions.\footnote{\label{iota}The $\iota$-operator means: `the only', $[\lambda y \ p]$ denotes: `(being) the object $y$ such that $p$', and $F=[\lambda y \ p]$ `the property $F$ constructed out of propositions $p$'.}

(Comp$_T$) may also be seen as a Comprehension Principle for \textit{structures} (that is, for \textit{models} of a theory $T$), taken to be the \textit{referents} of `true' theories. Thus, in practice, OT identifies theories with the \textit{structures} which satisfy them. $\mathsf{PA}$ is, for instance, the collection of all propositions true of $\mathsf{PA}$-\textit{structures}, $\mathsf{ZFC}$ the collection of all propositions true of $\mathsf{ZFC}$-structures, and so on. Also structures can be referred back to theories via (Comp$_T$). For instance, a model of $\mathsf{ZFC}$ which satisfies the Continuum Hypothesis (CH) will be interpreted as the abstract object corresponding to the theory $\mathsf{ZFC}$+CH. 

Overall, then, OT identifies models with \textit{particular theories} (and vice versa), and, as a consequence, it also provides a careful explanation of why $\mathsf{ZFC}$ is, for instance, different from $\mathsf{ZFC}$+CH: because the abstract object $\mathsf{ZFC}$+CH \textit{encodes} properties that are not encoded by the abstract object $\mathsf{ZFC}$. 

\medskip

OT also envisages the existence of a plenitude of first-order mathematical objects. OT's treatment of such a notion follows in a fully consequential manner from the theory's presuppositions. Let $\kappa$ be a mathematical object: the exact reference of $\kappa$ is fixed by some theory $T$. Thus, there can be no object $\kappa$ without an accompanying theory $T$ (which fixes its reference). So, $\kappa$ splits into a multitude of $\kappa_{T}$, one for each theory $T$. In formal terms: 

\[ \kappa_{T} \defeq \iota y (A!y \wedge (\forall F)(yF \leftrightarrow T \models F(\kappa_{T}))) \tag{Comp$_{obj}$} \]

\medskip

\noindent
(Comp), (Comp$_T$) and (Comp$_{obj}$) are sufficient for my purposes. In fact, only (Comp$_T$) shall prove indispensable. 

One last remark is in order. \textit{Encoding} also differs from \textit{exemplification}, insofar as abstract objects are `encoding-incomplete' whilst being `exemplification-complete'. This means that, say, a natural number $\kappa_{\mathsf{PA}}$, \textit{encodes} all and only those properties that $\mathsf{PA}$ asserts $\kappa$ to have. Therefore, it may be the case that $\kappa_{\mathsf{PA}}$ does not encode some property $F$, or its negation. On the other hand, for all $F$'s, either $\kappa_{\mathsf{PA}}$ exemplifies it or its negation. Thus, the encoding-related features of abstract objects are also able to painstakingly account for \textit{set-theoretic indeterminacy}. Reference to objects is fixed by theories, and these are, by their nature, `encoding-incomplete', so also first-order mathematical objects (\textit{sets}) will be encoding-incomplete. Hence, theories (of sets) are about inherently \textit{incomplete} objects, which are, however, liable to become `more complete'. A classic example is the real continuum (the power-set of $\omega$). The object $\mathcal{P}(\omega)_{\mathsf{ZFC}}$ is encoding-incomplete with respect to CH. But this doesn't mean that $\mathcal{P}(\omega)_T$ will, in general, be encoding-incomplete for any $T$. In fact, $\mathcal{P}(\omega)_{\mathsf{ZFC}+V=L}$ is encoding-complete with respect to CH. 

\section{Higher-Order Platonism Inside Object Theory}

HOP can be adequately interpreted (`embedded') inside OT. This is the present work's main claim, which I will now proceed to discuss. But first I should clarify what I really mean by this claim. 

My goal is not to suggest that Hamkins' multiverse is formally reducible (whatever that could mean) to a fragment of OT, nor do I posit that the Hamkinsian multiversist should think that his conception is just OT. What I actually think it that OT provides us with a detailed and coherent metaphysical account of HOP, which also enriches our understanding of Hamkins' multiverse.

In order to attain my goal, I will have to show that OT expresses the three main features of HOP, that is, (PLAT), (PERSP) and (PTP). This is a relatively trivial task as far as (PLAT) is concerned, but (PTP) and (PERSP) will need a more extended examination. 

\medskip

As far as (PLAT) is concerned: (Comp$_T$) envisages the existence of  a plenitude of higher-order abstract objects, that is, \textit{theories}. Moreover, as we have seen, (Comp$_T$) also entitles us to construing theories as structures satisfying them, and \textit{vice versa}. It is, in my view, a noticeable strength of OT, 
that it is also able to incorporate, and refer to, \textit{models} of set theory. To this end, as we have seen, one just has to pick out the relevant abstract objects corresponding to theories. 


(PTP) is also met by OT, since \textit{any} theory of sets represents a self-standing abstract object. Again, it is a defining feature of OT that theories are taken to be platonic objects. All of them are acceptable, that is, \textit{correct}, since OT's Comprehension Principles do not discriminate between `acceptable' and `unacceptable' \textit{entities}. Moreover, that any abstract object corresponding to a theory of sets $T$ reflects a `different concept of set', as envisaged by Hamkins, is expressed, in OT, by the fact that theories \textit{encode} different properties (in particular, they encode different theorems and truths about different structures). So, there is a clear and robust sense in which OT clearly asserts the existence of different `concepts of set'.   

Moreover, by using OT, now one may have at hand a way to formulate (PTP) which evades the objections in section \ref{HOP}: (PTP) now just reduces to OT's assertion that there exists a plenitude of theories, in particular, by (Comp$_T$) all \textit{conceivable} theories exist, so (PTP), in essence, just reduces to (Comp$_{T}$). In addition, one no longer needs to care about `consistent theories'. In fact, some theories sanctioned by (Comp$_{T}$) will, in fact, be inconsistent, which just means, by OT, that they are not exemplified, that is, that they do not have \textit{models}. I shall say more on this in section \ref{art}.


Finally, I proceed to address (PERSP), which, as I have said in section \ref{Intro}, is crucial to HOP's purposes. (PERSP) has been formulated as the view that the multiverse is always observed from the point of view of (is relative to) a specific model. In order to exemplify the position, I also mentioned Hamkins' own example of \textit{inner models of $V$}'s being best understood from the point of view of \textit{forcing extensions}. 

The other aspect of (PERSP) I have highlighted, and which is also to be reconsidered in light of OT, is that, within Hamkins' multiverse, one is supposed to constantly \textit{jump} from one model to another, and this is also construed (more broadly) by Hamkins as a `an enlargement of the theory/metatheory distinction'.    

In order to illustrate the first aspect, let us go back to the example of inner models inside forcing extensions of a provisional background universe $V$. One starts with constructing a forcing extension $V[G]$ of $V$ (for our, and Hamkins' purposes, this could just be any countable transitive model of $\mathsf{ZFC}$). The corresponding Hamkins multiverse, by (PERSP), is, thus, the collection of all models of $\mathsf{ZFC}$ accessible from $V[G]$. Among these, one finds inner models such as $L$, HOD, $L[\mathbb{R}]$, etc. Any of these will now reveal, inside $V[G]$, features which, in turn, are contingent on the features of $V[G]$ itself. 

Now, what kind of abstract objects and properties are, by OT, involved in observing, say, the constructible universe $L$ from the point of view of $V[G]$? The first object that needs consideration is $V[G]$ itself, that is, a model of $\mathsf{ZFC}$. Now, recall that, by (Comp$_T$), a model of $\mathsf{ZFC}$ is the \textit{same abstract object} as $\mathsf{ZFC}$. So, in particular, any property encoded by $V[G]$ will be \textit{encoded} by the abstract object $\mathsf{ZFC}$. But $V[G]$, presumably, encodes other properties that are not encoded by $\mathsf{ZFC}$, for instance, $\neg$CH (presumably, that was also the reason why we picked out $V[G]$ in the first place). In that case, we should, more correctly, refer $V[G]$ back to the abstract object (and theory) $\mathsf{ZFC}$+$\neg$CH.   

But now we need to express, in OT, that $\mathsf{ZFC}$+$\neg$CH does not only encode the features of $V[G]$, but also those of the $L$ \textit{inside} $V[G]$. So, our next object of scrutiny is the sentence:

\vspace{11pt}

\noindent
(\textbf{$\Phi$}). $V[G]$ is a model of $\mathsf{ZFC}$+$\neg$CH, and, moreover, $V[G] \models$ `$L$ has some set-theoretic property $\Psi$'.

\vspace{11pt}

\noindent
where $\Psi$ expresses some feature of $L$ as seen from the point of view of $V[G]$.\footnote{\cite{hamkins2012}'s example, on p. 418, is Cohen's construction of a $V[G]$ to get an inner model of $\neg$AC.} Now, $\Phi$ is a metatheoretic statement about models of $\mathsf{ZFC}$+$\neg$CH, so we cannot expect that it is really expressible in $\mathsf{ZFC}$+$\neg$CH. But, given the representability of $\mathsf{ZFC}$+$\neg$CH's metatheory in $\mathsf{ZFC}$+$\neg$CH itself, we may find, in the object theory, an equivalent statement, say, $\Phi_{\mathsf{ZFC+\neg CH}}$, that, somehow, corresponds to $\Phi$. Now, $\Phi_{\mathsf{ZFC+\neg CH}}$ is a \textit{bona fide} abstract object. In particular, by (Comp$_T$), $\Phi_{\mathsf{ZFC+\neg CH}}$ is a truth about some structure satisfying $\mathsf{ZFC}$+$\neg$CH itself, and, thus, can, finally, be viewed as being encoded by $\mathsf{ZFC}$+$\neg$CH. In the end, in OT, one has that the abstract object $\mathsf{ZFC}$+$\neg$CH encodes a property which expresses a `property of $L$ \textit{inside} a model $V[G]$ of $\mathsf{ZFC}$+$\neg$CH' (i.e., $\Psi$), and this is precisely what was required of me to establish.

In order to discuss the second aspect of (PERSP), i.e., the idea of `jumping' from one model to another, let us carry on with the example of $L$ inside $V[G]$. 

I seem to live in $L$ (in fact, the $L$ of $V[G]$), so my Hamkins multiverse consists of all models accessible from $L$ itself. On OT's conception, this just means that now I am carrying out my reasoning from the point of view of a different theory, that is, $\mathsf{ZFC}$+V=L, so, again, through (Comp$_T$), I am picking out a different abstract object, corresponding to the theory $\mathsf{ZFC}$+V=L. In other terms, the jump from $V[G]$ to $L$ can now be expressed as the shift from one abstract object to another, which is, in turn, tantamount to picking out (describing) \textit{distinct} abstract objects. If we wished to jump from there to another model $M$ (presumably of $\mathsf{ZFC}$), then we would just have to pick out the abstract object corresponding to the theory whose theorems are \textit{true} in $M$ and, of course, this process can be iterated as many times as one pleases.\footnote{I wish to stress that, on the one hand, `true in $M$' should be taken to be equivalent, through (Comp$_{T}$), to `true in (of) a theory $T$', wherein no specific model of $T$ is presupposed, and, on the other, that the satisfaction relation should always be seen as relative to $T$ itself (that is, should be more accurately symbolised as follows: `$\models_{T}$') and, as a consequence, to the models of $T$.} 

One further technical remark is necessary. Of course we cannot expect $\mathsf{ZFC}$+$\neg$CH (or $\mathsf{ZFC}$, for that matter) to literally prove that it has models (like $V[G]$) which satisfy $\neg$CH and further set-theoretic properties (of their $L$, for instance). The theory will, at most, prove, that some arbitrarily finite fragment of the theory has a \textit{countable} model with the mentioned properties, so the truth encoded by the abstract object corresponding to $\mathsf{ZFC}$+$\neg$CH will refer to such models. This is not ideal, but inevitable. On the other hand, if one used Hamkins' \textit{toy model perspective}, whereby all models are countable from the beginning, then one would get that the multiverse, and all relevant metamathematical facts about it, would be encoded only by the abstract object $\mathsf{ZFC}$.\footnote{For a discussion of the toy model perspective, see \cite{hamkins2012}, p. 436ff.} But this approach would be awkward, for then, in particular, the jump from one model to another would just be a jump from one countable model of $\mathsf{ZFC}$ to another countable model of $\mathsf{ZFC}$. Through keeping the principle of identifying a structure with the corresponding theory, we can, on the contrary, modulo the mentioned limitations, construe the jump from one model to another as picking out \textit{distinct} theories (abstract objects), something which seems to be a lot more faithful to the motivation behind Hamkins' multiverse. 

\medskip

In sum, OT is able to interpret (PLAT), (PERSP) and (PTP), which, overall, means that OT is able to interpret HOP. Let me empasise, again, that, in particular, OT is flexible enough to interpret the jumps from one model to another which are so characteristic of Hamkins' multiverse. The crux of OT's intepretation of (PERSP) consists in transforming metamathematical facts, on which the Hamkinsian multiversist thrives, into \textit{properties} encoded by abstract objects corresponding to the relevant theories. 
 
\section{Three Problems: Articulation, Skepticism, Practice}\label{SRI}

Let's take stock. I have concluded on a positive note about HOP's interpretability inside OT. I have also made it clear that OT should be taken to be a `companion' theory to, not a formalisation of, HOP. 

Now, HOP (and Hamkins' multiverse, for that matter) have sparked much controversy and raised concerns relating both to the underlying logic and the surrounding philosophy. In what follows, I will show that some of these concerns may be successfully addressed if one adopts the conception and tools afforded to us by OT.  

\subsection{Articulation}\label{art}

\cite{koellner2013} raises a problem of `articulation' for what he calls the `broad multiverse conception', that is, Hamkins' conception.\footnote{In fact, \cite{koellner2013}, pp. 4-5, distinguishes between the `broad multiverse conception', and the `relative broad multiverse conception', that is, the narrowing of the former to just one background theory $T$. For my purposes, it will just suffice to deal with the `broad multiverse conception'.} The problem may be explained as follows. 

The existence of a set-theoretic multiverse, that is, of all models of some theory of sets $T$, is secured through assuming the consistency of $T$. But, by G\"odel's Second Incompleteness Theorem, one cannot hope to prove Con($T$) in $T$, so, in order to get Con($T$), one needs to use a stronger theory as a background theory. Now, suppose one has ascended to a stronger theory $T'=T+Con(T)$ and can, now, articulate the multiverse of $T$. This would clearly be insufficient for the Hamkinsian multiversist, as he still cannot refer to models of $T'$, so one will have to ascend, again, to a stronger theory, $T''=T'+Con(T+Con(T))$ to articulate the multiverse of $T'$, and so on. Overall, this will result in an infinite regress through background theories of increasing consistency strength, none of which exhausts the `whole' Hamkinsian multiverse.\footnote{\cite{koellner2013}, pp. 7-8. A similar argument pointing out, this time, a `referential regress' in Hamkins' multiverse is in \cite{barton2016}, p. 202.}

The response to Koellner's argument is, in the light of OT, relatively straightforward. The HOP-ist does not need to articulate the `whole multiverse' in the way suggested by Koellner. In fact, by (PERSP), he doesn't even have a fixed concept of `set-theoretic multiverse'. It is true that he needs to start from some theory $T$, in order to build an initial multiverse, but, by (Comp$_T$), members of such multiverse will come to the fore automatically, since they are, so to speak, inbuilt features of the abstract object corresponding to the theory itself.

So, what the HOP-ist needs to do to articulate \textit{some} multiverse is just to pick out the abstract object corresponding to his initial theory $T$, and the properties this encodes. Among those properties, there will be the relevant metamathematical facts about the models of $T$. Afterwards, she can either pick out another theory (and abstract object) $T$ or, by (PERSP), she can progressively explore all other models accessible from the initial multiverse. 

Now, the HOP-ist's (PTP) may be construed as contravening the strategy above and, in fact, requiring the articulation of the \textit{whole} multiverse. But, as has been made clear in the previous section, inside OT, (PTP) just reduces to (Comp$_{T}$). Moreover, it should be stressed, once again, that the HOP-ist does not need to take into account the issue of the consistency of $T$, insofar as OT does not discriminate between consistent and inconsistent abstract objects. In particular, some abstract objects will encode inconsistent theories. However, by OT, this does not prevent such theories from existing: only, there are \textit{no} abstract objects which \textit{exemplify} them, so in particular there are \textit{no} structures which exemplify them.\footnote{As regards inconsistent objects, \cite{linsky-zalta1995}, p. 537, fn. 32, mentions the famous example of the `round square': by OT, there exists a unique abstract object encoding these two \textit{mutually inconsistent} properties, but there is \textit{no} object which exemplifies them. Note that, in principle, it is conceivable that there are \textit{no} abstract objects \textit{exemplifying} consistent theories, although this, presumably, would be contradicted by our experience of theories which we reasonably deem to be consistent, that is, such that these theories do, in fact, have \textit{models}.}

\medskip

\cite{koellner2013} also lays emphasis on what he thinks are further paradoxical features of Hamkins' multiverse relating to consistency.\footnote{Cf. \cite{koellner2013}, pp. 8ff.} Given any theory $T$ (of sufficient consistency strength), Hamkins' multiverse will also contain models of $T+\neg Con(T)$. Now, in order to secure the existence of those models in the multiverse, one has to presuppose that $T$ is, in fact, \textit{consistent}, so there is a tension between the `internal' and `external' claims of consistency of $T$. But, again, my interpretation of HOP gets rid of the distinction between internal and external. Access to abstract objects corresponding to theories is carried out through \textit{description}, so no \textit{external} presupposition needs to be made concerning the properties encoded by those objects, in particular, as we have seen, concerning their consistency. The existence of a (non-standard) model of $T+\neg Con(T)$ entirely lies `inside', and is encoded by, the abstract object $T+\neg Con(T)$, so the issue of the clash between the external and internal claims of consistency of such theory, practically, vanishes.

\subsection{Skepticism}

Does HOP instantiate a form of \textit{skepticism} about set theory as well as mathematical knowledge at large? 

In order to answer this question, one should, presumably, have at hand a clear notion of what would count as `skepticism', and the issue is too broad to even start examining it here. However, a few considerations on this issue may be inevitable, especially insofar as the charge has been levelled many times by different authors.\footnote{Cf. \cite{koellner2013}, pp. 9. \cite{button-walsh2018}, pp. 205ff., take  Hamkins' multiversism to be a contemporary version of Skolem's `model-theoretic skepticism', for which cf. \cite{skolem1922}.}    

One straightforward way of arguing that HOP is a form of skepticism is through raising the issue of \textit{referential indeterminacy} in connection with it. The argument, broadly, runs as follows: set-theorists seem to have \textit{determinate} (\textit{uni-que}) referents in mind for set-theoretic objects, whereas HOP subscribes to (Comp$_{obj}$), which implies, among other things, that referents of set-theoretic objects are \textit{non-unique} in an essential way; hence HOP cannot represent an accurate interpretation of the set-theoretic discourse, let alone a foundation for it. 

One possible response is that this argument may be thriving on a (subtle) confusion between `relativism' and `indeterminacy'. Granted, HOP is relativistic, since properties of set-theoretic objects are deemed by it to be contingent on theories (and structures modelling them). Yet, it would not be correct to view HOP as supporting the referential indeterminacy of the set-theoretic discourse. HOP (and OT, for that matter), do not hold that the reference of, say, 1 is \textit{indeterminate}; rather, they hold that it is \textit{plural}. According to HOP-ists, and OT-ists as well, 1 is an incomplete denotation of the object $1_{T}$, that is, of 1 as defined within some theory $T$. Now, they may agree that, when mathematicians (and set-theorists) refer to 1, they refer, most of the times, to 1$_{\mathsf{PA}}$, but, will, as is clear, deny that this is inevitable. For instance, 1$_{\mathsf{ZF}}$ is a sensible alternative that comes with a lot of new theory.\footnote{However, the distinction between, say, $\kappa_{\mathsf{T_1}}$, the object $\kappa$ as seen from the point of view of a theory $T_1$, and $\kappa_{\mathsf{T_2}}$, the object $\kappa$ as seen from the point of view of another theory $T_2$, should, sometimes, be carefully delimited, if (PERSP) is to be implemented satisfactorily. For instance, one may want to say that $\kappa_{\mathsf{T_1}}$ is countable, and $\kappa_{\mathsf{T_2}}$  uncountable. Then, for the purpose of comparing them, it would be desirable to have that $\omega_{\mathsf{T_1}}$ is the same as $\omega_{\mathsf{T_2}}$. However, it is not necessary to further posit that `$\cong_{T_1}$' is the same as `$\cong_{T_2}$', as the use of the `$\cong$' relationship is clearly bounded by the domains of quantification of the respective theories.}   

One could insist that the interpretation of the arithmetical discourse is so plain and evident, that there is no reason to assume that the referent of a natural number is contingent on a specific theory, but I think that this only begs the question.

As a consequence, HOP-ists' appeal to the relativity of the set-theoretic discourse should be more correctly construed as expressing the belief in the \textit{multiple realisability}, not in the \textit{indeterminacy}, of the discourse itself.  Hamkins has stressed this point many times:

\begin{quote}

[...] let me remark that this multiverse vision, in contrast to the universe view with which we began this article, fosters an attitude that what set theory is about is the exploration of the extensive range of set-theoretic possibilities. (\cite{hamkins2012}, p. 440)
    
\end{quote}

\noindent
As regards the issue of whether such a view instantiates any form of skepticism, he explains that:

\begin{quote}

..to my way of thinking, this label is misapplied, for the multiverse
position is not especially or necessarily skeptical about set-theoretic
realism. We do not describe a geometer who works freely sometimes in Euclidean and sometimes in non-Euclidean geometry as a geometry skeptic. (\cite{hamkins2020}, p. 295-6)

\end{quote}

\noindent
Now, it seems to me that my interpretation of HOP \textit{within} OT adds crucial weight to this argument. Believers in the existence of all conceivable abstract objects do not deny the \textit{existence} (or \textit{knowability}) of (mathematical) reality. They think, instead, that mathematical reality consists of a very rich realm of realit\textit{ies}, from which one can pick one's preferred mathematical concepts. Inside OT, this idea prominently features as the claim that mathematical objects are \textit{non-spatiotemporal} and \textit{sparse}.\footnote{On this, cf., \cite{linsky-zalta1995}, especially section VI.} 

The accusation frequently levelled against \textit{plenitudinous platonism} of being a brand of formalism, that is, of `anti-realism', is not new in the debate, in any case, and discussing it in full is besides the scope of this paper.\footnote{For instance, \cite{field2001}, pp. 334ff. takes Balaguer's full-blooded platonism (FBP) to be of an `anti-objectivistic' character, and \cite{potter2004}, p. 11, even denies that FBP has anything to do with platonism. \cite{linnebo2018}, p. 21, on the other hand, views such conceptions as FBP and OT as belonging `somewhere in the territory between anti-nominalism and full-fledged platonism'.} Presumably, whether one may have a coherent concept of `multiple (mathematical) realities' also depends on what one takes reality to consist, and on how the latter (causally or non-causally) interacts with human minds. Presumably on such grounds, which I cannot examine here, some people will keep insisting that there is a \textit{unique} mathematical reality and, therefore, that plenitudinous platonism (as well as HOP) have a skeptical attitude about set-theoretic reality.\footnote{Cf., for instance, \cite{martin2001}.} But this is a claim we ought not to light-heartedly grant.  

\subsection{Practice}\label{meta}

Our theorising, so far, has been `metaphysically thick'. By HOP, in order to understand what set-theorists really do, one does not only need to understand set theory mathematically, but also understand some accompanying, explanatory theory of set-theoretic objects, theories and structures like OT. 

Thus, one further objection against HOP would be that one doesn't actually need to do this. That is, one could say that, whatever the merits of the position, regardless of whether it's coherent or not, such an extreme form of set-theoretic realism really has no impact on set-theorists' work. Moreover, the objector would continue, the multiversist seems to want to make things more complicated than they are, since taking onboard a complex, platonistically construed, multiverse ontology affords no potential gains in foundational terms. 

Naturalists, such as Maddy, have expressed this point of view. In particular, Maddy has stressed that Universism \textit{and} the theory $\mathsf{ZFC}$+Large Cardinals are more than sufficient for all our mathematical and philosophical purposes, where the role of philosophy is taken, by Maddy, to consist in that of Second Philosophy, that is, of a philosophy which attends to the `details of practice' and does not aim to reinterpret or influence in any way the practitioners' work.\footnote{For the articulation of Second Philosophy, cf. \cite{maddy2007}.} 

For instance, \cite{maddy2017} denounces all metaphysical incursions into mathematical work, by declaring that:

\begin{quote}

[t]he metaphysics of abstracta or meanings or concepts are all really beside the point. The fundamental challenge these multiverse positions raise for the universe advocate is this: are there good reasons to pursue a single, preferred theory of sets that’s as decisive as possible, or are there not? (\cite{maddy2017}, p. 316)

\end{quote}

\noindent
We know Maddy's answer to her own question: the reasons to pursue just one theory, $\mathsf{ZFC}$+Large Cardinals and `single-$V$', are more cogent than the reasons to pursue a multiverse theory, since $\mathsf{ZFC}$+Large Cardinals and `single-$V$' are taken to best express the foundational goals of set theory. 

Now, although the universist orthodoxy, so to speak, is far from being overthrown, several works have appeared which put pressure on this conception. For instance, \cite{ternullo2019} suggests that the goals mentioned by Maddy might also be attained by a multiverse theory, and that consideration of further fundamental set-theoretic goals, such as studying the relationships between models (theories), might, in fact, lead one to prefer a multiverse over a single-universe theory.\footnote{Cf. \cite{ternullo2019}, pp. 57ff.} Moreover, `multiverse maths' has delivered many promising mathematical programmes, and many more might appear in the next future.\footnote{An immediate example that comes to mind is set-theoretic `geology', for which see \cite{fuchs-hamkins-reitz2015}, now thought to be central to many set-theoretic undertakings.}

More recently, \cite{antos2022} has pointed out that also a Second Philosopher could see models of set theory as being part of set theory's `fundamental entities'. This is because:

\begin{quote}

[i]ntroducing models as fundamental entities is part of set-theoretic methodology and not just a heuristic aid. Being about methodology, this claim then becomes eligible to inform further philosophical questions [...]. (\cite{antos2022}, p. 6) 
    
\end{quote}

\noindent
Of course, at this point, one could observe that one thing is to hold that the set-theoretic multiverse (or some conception thereof) might fit with a purely naturalistic methodology, and another is to advocate a complex, practically dispensable, metaphysical theory such as OT. 

In other terms, while the naturalist could, in principle, be open to embrace \cite{mostowski1967}'s formalistic pluralism about `set theory', as expressed in this quote:

\begin{quote}

Probably we shall have in the future essentially different intuitive notions of sets just as we have different notions of space, and will base our discussions of sets on axioms which correspond to the kind of sets which we want to study. (\cite{mostowski1967}, p. 94)

\end{quote}

\noindent
she could not possibly subscribe to the metaphysical theorising inherent in HOP. At which point, the HOP-ist could reply that HOP is consistent with the kind of naturalism that includes both Mostowski's pluralism and Antos' suggestion that models are indispensable mathematical entities. As a consequence, the anti-platonist naturalist will now turn to launch a full-scale attack against HOP's platonistic doctrines, and the HOP-ist will respond, as Hamkins does, that platonism is `natural', insofar as the \textit{experience} of models of set theory is genuine, and cannot be explained away.\footnote{\cite{hamkins2012}, especially p. 418. This could be expressed more intelligibly as follows: there is no \textit{alternative} explanation of the existence of a plurality of models which is as good as the HOP-ist explanation; in particular, universist explanations, which hold that models are just convenient fictions (inside a single $V$), won't do, as the experience of models is so mathematically substantial that considering the existence of models a matter of fact becomes, by the Maddian naturalist's own lights, inevitable.} Then, the naturalist will respond that the HOP-ist's notion of `experience' is not hers, and so on. This debate will, most likely, end in a stalemate.

It might, in any case, be too early to assess whether HOP's foundation of set theory is going to fully convince the naturalistically-minded philosophers and set-theorists (incidentally, this might also depend on the status of Universism, as well as of Universism-inspired mathematical programmes, in the coming years).

What is certain is that, through OT, HOP seems to have noticeable advantages over other conceptions: correct or not, it provides an \textit{explanation} of set-theoretic \textit{undecidability}, of the fact, that is, that many set-theoretic statements cannot be decided by the current set-theoretic axioms, and of why undecidability will persist, by precisely bringing out the  \textit{incomplete}, in particular, \textit{encoding-incomplete}, character of set-theoretic concepts and theories qua OT's abstract objects, i.e., the fact that these objects do not have enough mathematical `content' to decide lots of set-theoretic statements. It might not be too much for a naturalist, but having at hand some robust explanation of set-theoretic undecidability will sound appealing to many.

Finally, \cite{clarke-doane2022} sharply points out that mathematical pluralism (one dare say, by now, HOP) is also able to account for the \textit{reliablity} of our mathematical knowledge. The `reliability challenge' is the challenge to the view that we have \textit{true} mathematical knowledge. Mathematical pluralism, in a sense, dissolves the problem: it simply isn't possible, in the light of mathematical pluralism, that `we could have had \textit{false} mathematical knowledge' since pluralists commit themselves to the idea that any consistent collection of set-theoretic propositions (theory) instantiates \textit{some} truths.\footnote{Cf. \cite{clarke-doane2022}, especially Ch. 3.}

\section{Concluding Remarks}

I have shown that there is a way to make sense of `higher-order realism' in the context of Hamkins' multiverse, by (informally) embedding what I have called HOP into a general theory of abstract objects, OT. Such an embedding does not require of one to see HOP as a `fragment' of OT, only to acknowledge that HOP's platonistic conception can be adequately expressed by OT. 

In particular, the Hamkinsian multiversist's typical statements that \textit{models} are set-theorists' main objects of investigation, that we have a genuine experience of them, that set-theorists jump from one universe to another, become fully intelligible in the context of OT. Moreover, as we have seen, OT comes with an explanation of set-theoretic undecidability, which, in my view, makes the `multiverse undertaking' even more persuasive. Of course, one can resist the force of the explanation, but, certainly, the latter stands Hamkins' multiversism in good stead. 

Finally, I have also shown that, through OT's reinterpretation of HOP, the Hamkinsian multiversist may also be able to respond to issues concerning the articulation of the multiverse, and reject as, at least in part, inappropriate, the label of `skepticism' for his conception.

OT presents itself as an elegant (and economic) solution to the interpretative problems posed by (PLAT), (PERSP) and (PTP). The basic objects are taken to be theories, and these are reduced to the structures which satisfy them. So the complex metatheory of the multiverse which is needed to express salient facts about the models of set theory can be incorporated, through the very flexible and apt notion of \textit{encoding}, into the theories themselves.

As hinted at at the beginning of section \ref{HOP} and in fn. \ref{plen}, there are other accounts of higher-order platonism out there in the market, and it cannot be ruled out that further, equally suitable, platonistic accounts of models of set theory, of multiversism at large, will be found.\footnote{For instance, although not directly addressing set-theoretic pluralism, \cite{horsten2019}'s metaphysical account of `arbitrary objects' seems very promising in this respect. For a comparison between  Horsten's own approach and OT, see section 8.3 of that work.} On the contrary, it would be entirely legitimate to wonder:

\begin{question}

In which ways could the study of the metaphysics of HOP be further developed, and with what foundational implications? 
    
\end{question}

\noindent
Once shown that it is concretely possible to bring the metaphysical features of the multiverse to bear on the foundations of set theory, where does all this leave us? Does such incursion into metaphysics imply some major revision of our accepted and practiced methodologies in set theory? 

My answer has been, and still is, a `no', especially insofar as HOP does not directly bear on the mathematical methodologies in any way. But I think that the potential usefulness of `auxiliary conceptions' in the context of the foundations of mathematics should not be denied straight away. This is because, as Hamkins himself asserts: 

\begin{quote}
    ..our foundational ideas in the mathematics and philosophy of set theory outstrip our formalism, mathematical issues become philosophical, and set theory increasingly finds itself in need of philosophical assistance. (\cite{hamkins2012}, p. 436)
\end{quote}

\noindent
although, admittedly, the dividing line between `philosophical assistance' and `philosophical reinterpretation of mathematics' might, potentially, be very hard to draw. 

This last remark very aptly leads me to pose one further question: 

\begin{question}

To what extent should extra-set-theoretic, metaphysical, in particular, resources be taken to bear on set-theoretic work? 
    
\end{question}

\noindent
At this stage, I could not say anything on this issue which would not look like as just a collection of superficial remarks. I will just limit myself to observing that the present work may, potentially, contribute to the view that, once the extra-set-theoretic resources are properly formalised and adequately interpreted, they might really become part of (even explain) the underlying logic and concepts of set theory in a meaningful way, and not just play the role of a `useful heuristic', in the Maddian sense. But clearly it would still be a long way to go to make a fully convincing case for this claim.

\pagebreak

\bibliography{Bib1}

\begin{thebibliography}{}

\bibitem[Antos, 2022]{antos2022}
Antos, C. (2022).
\newblock Models as {F}undamental {E}ntities in {S}et {T}heory: {A}
  {N}aturalistic and {P}ractice‑based {A}pproach.
\newblock {\em Erkenntnis}, pages 1--28.

\bibitem[Antos et~al., 2015]{afht2015}
Antos, C., Friedman, S.-D., Honzik, R., and Ternullo, C. (2015).
\newblock Multiverse {C}onceptions in {S}et {T}heory.
\newblock {\em Synthese}, 192(8):2463--2488.

\bibitem[Bagaria, 2005]{bagaria2005}
Bagaria, J. (2005).
\newblock Natural {A}xioms of {S}et {T}heory and the {C}ontinuum {P}roblem.
\newblock In H\'ajek, P., Vald{\'e}s-Villanueva, L., and Westert{\aa}hl, D.,
  editors, {\em Logic, {M}ethodology and {P}hilosophy of {S}cience:
  {P}roceedings of the {T}welfth {I}nternational {C}ongress}, pages 43--64.
  King's College Publications, London.

\bibitem[Balaguer, 1995]{balaguer1995}
Balaguer, M. (1995).
\newblock A {P}latonist {E}pistemology.
\newblock {\em Synthese}, 103:303--25.

\bibitem[Barton, 2016]{barton2016}
Barton, N. (2016).
\newblock Multiversism and concepts of set: {H}ow much relativism is
  acceptable?
\newblock In Boccuni, F. and Sereni, A., editors, {\em Objectivity, {R}ealism
  and {P}roof}, pages 189--209. Springer, Berlin.

\bibitem[Bernays, 1935]{bernays1935}
Bernays, P. (1935).
\newblock Sur le {P}latonisme dans les {M}ath{\'e}matiques.
\newblock {\em L'Enseignement Math{\'e}matique}, 34:52--69.

\bibitem[Blechschmidt, 2022]{blechschmidt2022}
Blechschmidt, I. (2022).
\newblock Exploring {M}athematical {O}bjects from {C}ustom-{T}ailored
  {M}athematical {U}niverses.
\newblock In Oliveri, G., Ternullo, C., and Boscolo, S., editors, {\em Objects,
  {S}tructures and {L}ogics. {F}ilMat {S}tudies in the {P}hilosophy of
  {M}athematics}, Boston Studies in Philosophy of Science, pages 63--95.
  Springer, Cham.

\bibitem[Button and Walsh, 2018]{button-walsh2018}
Button, T. and Walsh, S. (2018).
\newblock {\em Philosophy and {M}odel {T}heory}.
\newblock Oxford University Press, Oxford.

\bibitem[Clarke-Doane, 2022]{clarke-doane2022}
Clarke-Doane, J. (2022).
\newblock {\em Mathematics and {M}etaphilosophy}.
\newblock Cambridge University Press, Cambridge.

\bibitem[Field, 2001]{field2001}
Field, H. (2001).
\newblock {\em Truth and {A}bsence of {F}act}.
\newblock Oxford University Press, Oxford.

\bibitem[Fuchs et~al., 2015]{fuchs-hamkins-reitz2015}
Fuchs, G., Hamkins, J., and Reitz, J. (2015).
\newblock Set-{T}heoretic {G}eology.
\newblock {\em Annals of {P}ure and {A}pplied {L}ogic}, 166(4):464--501.

\bibitem[G{\"o}del, 1947]{godel1947}
G{\"o}del, K. (1947).
\newblock What is {C}antor's {C}ontinuum {P}roblem?
\newblock {\em American Mathematical Monthly}, 54:515--525.

\bibitem[G{\"o}del, 1951]{godel1951}
G{\"o}del, K. (1951).
\newblock Some {B}asic {T}heorems on the {F}oundations of {M}athematics and its
  {A}pplications.
\newblock Gibbs Lecture.

\bibitem[Hamkins, 2012]{hamkins2012}
Hamkins, J.~D. (2012).
\newblock The {S}et-{T}heoretic {M}ultiverse.
\newblock {\em Review of Symbolic Logic}, 5(3):416--449.

\bibitem[Hamkins, 2020]{hamkins2020}
Hamkins, J.~D. (2020).
\newblock {\em Lectures on the {P}hilosophy of {M}athematics}.
\newblock MIT Press, Cambridge (MA).

\bibitem[Horsten, 2019]{horsten2019}
Horsten, L. (2019).
\newblock {\em The {M}etaphysics and {M}athematics of {A}rbitrary {O}bjects}.
\newblock Cambridge University Press, Cambridge.

\bibitem[Isaacson, 2011]{isaacson2011}
Isaacson, D. (2011).
\newblock The {R}eality of {M}athematics and the {C}ase of {S}et {T}heory.
\newblock In Nov\'ak, Z. and Simonyi, A., editors, {\em Truth, {R}eference and
  {R}ealism}, pages 1--75.

\bibitem[Jech, 2003]{jech2003}
Jech, T. (2003).
\newblock {\em {S}et {T}heory}.
\newblock Springer, Berlin.

\bibitem[Koellner, 2013]{koellner2013}
Koellner, P. (2013).
\newblock Hamkins on the {M}ultiverse.
\newblock Pre-print.

\bibitem[Kunen, 2011]{kunen2011}
Kunen, K. (2011).
\newblock {\em {S}et {T}heory. {A}n {I}ntroduction to {I}ndependence {P}roofs}.
\newblock College Publications, London.

\bibitem[Lewis, 1986]{lewis1986}
Lewis, D. (1986).
\newblock {\em On the {P}lurality of {W}orlds}.
\newblock Blackwell, Oxford.

\bibitem[Linnebo, 2017]{linnebo2017}
Linnebo, {\O}. (2017).
\newblock {\em Philosophy of {M}athematics}.
\newblock Princeton University Press, Princeton.

\bibitem[Linnebo, 2018]{linnebo2018}
Linnebo, {\O}. (2018).
\newblock Platonism in the {P}hilosophy of {M}athematics.

\bibitem[Linsky and Zalta, 1995]{linsky-zalta1995}
Linsky, B. and Zalta, E. (1995).
\newblock Naturalized {P}latonism versus {P}latonized {N}aturalism.
\newblock {\em Journal of Philosophy}, XCII(10):525--555.

\bibitem[Linsky and Zalta, 2006]{linsky-zalta2006}
Linsky, B. and Zalta, E. (2006).
\newblock What is {N}eologicism?
\newblock {\em The Bulletin of Symbolic Logic}, 12(1):60--99.

\bibitem[Maddy, 2007]{maddy2007}
Maddy, P. (2007).
\newblock {\em Second {P}hilosophy: a {N}aturalistic {M}ethod}.
\newblock Oxford University Press, Oxford.

\bibitem[Maddy, 2017]{maddy2017}
Maddy, P. (2017).
\newblock Set-{T}heoretic {F}oundations.
\newblock In Caicedo, A., Cummings, J., Koellner, P., and Larson, P.~B.,
  editors, {\em Foundations of {M}athematics. {E}ssays in {H}onor of {W}.
  {H}ugh {W}oodin's 60th {B}irthday}, Contemporary Mathematics, 690, pages
  289--322. American Mathematical Society, Providence (Rhode Island).

\bibitem[Martin, 2001]{martin2001}
Martin, D. (2001).
\newblock Multiple {U}niverses of {S}ets and {I}ndeterminate {T}ruth {V}alues.
\newblock {\em Topoi}, 20:5--16.

\bibitem[Mostowski, 1967]{mostowski1967}
Mostowski, A. (1967).
\newblock Recent {R}esults in {S}et {T}heory.
\newblock In Lakatos, I., editor, {\em Problems in the {P}hilosophy of
  {M}athematics. {P}roceedings of the {I}nternational {C}olloquium in the
  {P}hilosophy of {S}cience, {L}ondon, 1965}, volume~47, pages 82--108.
  Elsevier, Amsterdam.

\bibitem[Nodelman and Zalta, 2014]{nodelman-zalta2014}
Nodelman, U. and Zalta, E. (2014).
\newblock Foundations for {M}athematical {S}tructuralism.
\newblock {\em Mind}, 123(489):39--78.

\bibitem[Panza and Sereni, 2013]{panza-sereni2013}
Panza, M. and Sereni, A. (2013).
\newblock {\em Plato's {P}roblem. {A}n {I}ntroduction to {M}athematical
  {P}latonism}.
\newblock Palgrave-MacMillan, London.

\bibitem[Parsons, 1995]{parsons1995}
Parsons, C. (1995).
\newblock Platonism and {M}athematical {I}ntuition in {K}urt {G}{\"o}del's
  {T}hought.
\newblock {\em Bulletin of Symbolic Logic}, 1(1):44--74.

\bibitem[Potter, 2004]{potter2004}
Potter, M. (2004).
\newblock {\em {S}et {T}heory and its {P}hilosophy}.
\newblock Oxford University Press, Oxford.

\bibitem[Putnam, 1980]{putnam1980}
Putnam, H. (1980).
\newblock Models and {R}eality.
\newblock {\em Journal of Symbolic Logic}, 45(3):464--82.

\bibitem[Shapiro, 1997]{shapiro1997}
Shapiro, S. (1997).
\newblock {\em Philosophy of {M}athematics. {S}tructure and {O}ntology}.
\newblock Oxford University Press, Oxford.

\bibitem[Shelah, 2003]{shelah2003}
Shelah, S. (2003).
\newblock Logical {D}reams.
\newblock {\em Bulletin of the American Mathematical Society}, 40(2):203--228.

\bibitem[Skolem, 1967]{skolem1922}
Skolem, T. (1967).
\newblock Some remarks on axiomatized set theory.
\newblock In van Heijenoort, J., editor, {\em From {F}rege to {G}\"odel. A
  {S}ource {B}ook in {M}athematical {L}ogic, 1879-1931}, pages 290--301.
  Harvard University Press, Cambridge (MA).

\bibitem[Ternullo, 2019]{ternullo2019}
Ternullo, C. (2019).
\newblock Maddy on the {M}ultiverse.
\newblock In Centrone, S., Kant, D., and Sarikaya, D., editors, {\em
  Reflections on the {F}oundations of {M}athematics: {U}nivalent {F}oundations,
  {S}et {T}heory, and {G}eneral {T}houghts}, pages 43--78. Springer Verlag,
  Berlin.

\bibitem[Woodin, 2001]{woodin2001}
Woodin, W.~H. (2001).
\newblock The {C}ontinuum {H}ypothesis.
\newblock {\em Notices of the American Mathematical Society}, Part 1: 48, 6, p.
  567--76; Part 2: 48, 7, p. 681--90.

\bibitem[Woodin, 2011]{woodin2011c}
Woodin, W.~H. (2011).
\newblock The {C}ontinuum {H}ypothesis, the {G}eneric-{M}ultiverse of {S}ets,
  and the ${\Omega}$-{C}onjecture.
\newblock In Kennedy, J. and Kossak, R., editors, {\em Set {T}heory,
  {A}rithmetic, and {F}oundations of {M}athematics: {T}heorems,
  {P}hilosophies}, pages 13--42. Cambridge University Press, Cambridge.

\bibitem[Zalta, 1983]{zalta1983}
Zalta, E. (1983).
\newblock {\em Abstract {O}bjects: an {I}ntroduction to {A}xiomatic
  {M}etaphysics}.
\newblock Reidel, Dordrecht.

\bibitem[Zalta, 2000]{zalta2000}
Zalta, E. (2000).
\newblock Neo-logicism? {A}n {O}ntological {R}eduction of {M}athematics to
  {M}etaphysics.
\newblock {\em Erkenntnis}, 53:219--65.

\end{thebibliography}

\bibliographystyle{apalike}

\pagebreak

\tableofcontents 

\end{document}